# Virasoro algebra and Sugawara constructions formally in Lean

Kalle Kytölä

Aalto University, Department of Mathematics and Systems Analysis
kalle.kytola@aalto.fi

**Abstract**

We formalize in Lean certain calculational proofs about infinite-dimensional Lie algebras. Specifically, we construct the Virasoro algebra as a central extension of the Witt algebra associated with a nontrivial 2-cocycle, and we construct representations of the Virasoro algebra by Sugawara constructions.

## 1 Introduction

This report describes a project to formalize the basics of the Virasoro algebra and Sugawara constructions in the Lean theorem prover. The report is intended for mathematicians and formalizers. The main assumed mathematical background is familiarity with the definition of a Lie algebra, which also served as the essential starting point in Lean's Mathlib for this project. While the text aims to be readable with a minimal knowledge of Lean syntax, some familiarity with formalization of mathematics may be necessary to appreciate the implementation and design considerations.

A mathematician could use this report to get a feeling for what the formalization of some standard mathematical topics entails. Together with the accompanying online blueprint, the report should also serve as a fairly detailed account of the calculation of the 2-cohomology of the Witt algebra, the definition of the Virasoro algebra, and the calculational proofs of variants of the Sugawara construction.

To the working formalizer, this experiment mainly offers the reassurance that Mathlib serves well as a foundation on which to build more mathematics (gradually approaching research level), and that the formalization community's wisdom largely continues to apply. We discuss some formalization issues relevant to calculational topics in particular, and a few specific formalization design experiments.

Table 1 summarizes the main objects and results formalized in the project, and indicates where in this report they are discussed. The project repository

https://github.com/kkytola/VirasoroProject

hosts the Lean code, in total approximately 4000 lines currently. The web page

https://kkytola.github.io/VirasoroProject/

contains the documentation and a blueprint, including a dependency graph, which can serve as a fairly detailed (informal) guide to the project's contents. For definiteness, this report discusses progress at the time of commit 29e118b; the project will be updated and expanded later.

### 1.1 Context for the mathematics

The Virasoro algebra [Vir70, Vir08] is an infinite-dimensional Lie algebra which plays a central role in two-dimensional conformal field theory (CFT) and string theory [BPZ84]. Conformal transformations in two dimensions are infinitesimally described by the Witt algebra of holomorphic vector fields, but since symmetries in quantum theories are represented projectively, central extensions of literal symmetry algebras enter. The Virasoro algebra is the unique such nontrivial central extension of the Witt algebra [GF69].

A number of fundamentally important conformal field theories possess more primitive symmetries and associated "currents" [Sug68], and the Virasoro algebra actions can be constructed from these primitive objects via so-called Sugawara constructions [GO86]. The details of Sugawara constructions vary, but the common theme is that the generators of the Virasoro algebra can be expressed as formally infinite quadratic expressions in the generators of the current algebra (reflecting the fact that the energy-momentum tensor is quadratic in the currents). Examples of such theories include the (massless) free boson [Vir70, FK81], free fermions [Ram71, NS71], symplectic fermions [GK99], Wess–Zumino–Witten models [FK81, GO85], and and their coset theories [GKO86].

### 1.2 Lean, Mathlib, and formalization

The formalization is done in *Lean 4* [MU21], a dependently typed language and a proof assistant. We build on Lean's formalized mathematical library *Mathlib* [mat20]. Mathlib is a maintained open-source library of very high standards, which aims to serve as a foundation for essentially all areas of mathematics, to have definitions and results in a very general form, and to follow unified conventions, so that different parts of the library work seamlessly together and (eventually) enable the formalization of mathematical research. Mathlib's current scope is a good coverage of undergraduate mathematics (at a level of generality and abstraction that far exceeds undergraduate teaching), a substantial amount of graduate level mathematics, and even some essentially research level topics. The high standards of Mathlib come at a cost of inertia in the contributions process. The results here have not yet been pull-requested to Mathlib.

| Topic | Reference | Lean declarations and links to documentation |
|---|---|---|
| central extensions of Lie algebras | Def 3.1 | `LieAlgebra.IsCentralExtension` <br> `LieTwoCocycle.CentralExtension` |
| Lie algebra cohomology in degree 2 | Def 4.3 | `LieTwoCohomology` |
| Virasoro algebra | Def 8.1 | `VirasoroAlgebra` |
| Heisenberg algebra | Sec 9 | `HeisenbergAlgebra` |
| Verma modules | Def 12.2 | `VermaModule` <br> `IsVermaModule` |
| calculation of the 2nd cohomology group of the Witt algebra | Thm 6.2 | `WittAlgebra.rank_lieTwoCohomology_eq_one` |
| Sugawara construction for Heisenberg algebra | Thm 10.1 | `sugawaraRepresentation_of_representation_heisenbergAlgebra` |
| module map from Virasoro Verma module to charged Fock space | Thm 12.7 | `ChargedFockSpace.virasoroVermaToChargedFockSpace` |

**Table 1.** Main objects and results formalized and discussed in this report.

The current formalization project was not intended as systematic library building, but rather as an experiment to use the current Mathlib as a basis for the formalization of a specific topic, which in turn is foundational to some active areas of research. We see such fleshing out of particular topics as meaningful tests of a general library of mathematics.
We include Lean code excerpts throughout to exemplify the formalization and to make specific remarks about the design choices. The excerpts are essentially directly from the repository's current codebase (and we include links to it), with some liberties taken for exposition.[1] For exemplified definitions (`def` or `abbrev`), we sometimes display the formal construction, whereas for results (`lemma` or `theorem`), we usually only show the formal statement of the result and omit the formal proof. The full formalizations are available in the public repository. For a reader not familiar with Lean, a quick introduction to some essential syntax is provided in Figure 1.

### 1.3 Related work

We are not aware of earlier formalizations in any proof assistants of any substantial part of using Lie algebra cohomology to classify central extensions, calculating cohomology groups of specific (infinite-dimensional) Lie algebras explicitly, constructing the Virasoro algebra and Heisenberg algebra as central extensions, or constructing representations of the Virasoro algebra by Sugawara-type constructions.

The project directly builds on the formalization of Lie algebras in Mathlib contributed primarily by Oliver Nash [Nas22]. Nash's work included constructions of specific Lie algebras, including all finite-dimensional complex semisimple Lie algebras.

Group cohomology has been developed in Lean and contributed to Mathlib [Liv23], and the formalization of cohomology theories generally is progressing [BBCM23]. In [Liv23], group cohomology was built in general, and then supplemented with specific API for cohomology in low degrees, which has direct applications. The application of group cohomology in degree two to group extensions is not yet in Mathlib. The present project, by contrast, essentially starts by constructing interesting central extensions of Lie algebras corresponding to Lie algebra cohomology in degree two. We only wrote the low-degree cohomology and API relevant for that purpose; cohomology in general degrees or with coefficients in general modules was not done.

Essentially concurrently with the present work, pull request #20206 by Scott Carnahan on Lie algebra extensions was merged to Mathlib. The definition is almost directly compatible with the current project's one, so reconciling the small overlap should be easy.

Some ambitious formalization work of Carnahan [Car25] builds towards the theory of vertex operator algebras (VOA). These are a far more advanced topic related to breakthroughs in group theory [FLM89], and also a starting point for a mathematical definition of conformal field theory [Gui23]. VOAs include the Virasoro algebra as a component, often via variants of Sugawara constructions. We hope that the current project finds synergies with the formalization of VOAs.

## 2 Lie algebras

We assume familiarity with Lie algebras, but include the definition for completeness, for comments about the Lean implementation, and for specific later references.

**Definition 2.1.** A **Lie algebra** over a field $\mathbb{F}$ is a vector space $\mathfrak{g}$ over $\mathbb{F}$, equipped with a bilinear **Lie bracket**

---

[1] Besides deliberate exposition choices, some Unicode symbols had to be typeset in clumsier similar alternatives, which are not literally correct Lean.

Virasoro algebra and Sugawara constructions formally in Lean

```
-- Let `R` be a (semi)ring and let `V` be an `R`-module.
variable {R : Type*} [Semiring R] {V : Type*} [AddCommGroup V] [Module R V]
                    standing hypotheses/context

/-- Commutator `[A,B] := AB-BA` of two linear operators. -/
def LinearMap.commutator (A B : V →₁[R] V) : V →₁[R] V := A * B - B * A
                          context/hypotheses     definition type   definition content

/-- `[A,B] = -[B,A]` -/
lemma LinearMap.commutator_comm (A B : V →₁[R] V) : A.commutator B = - B.commutator A := by
  simp [LinearMap.commutator]
                                    context/hypotheses          lemma statement
  proof of lemma
```

**Figure 1.** Lean syntax exemplified: The keyword `variable` introduces variables and standing assumptions. The keywords `def` and `abbrev` make new definitions. The two functionally synonymous keywords `lemma` and `theorem` are for results and their proofs. Definitions, lemmas, and theorems are named (e.g. `LinearMap.commutator` above). For theorems (resp. definitions), after the name come the hypotheses (resp., the context in which a definition is made), then after a colon `:` the theorem statement (resp. the type of defined object), and finally after `:=` the proof (resp. the actual implementation of the definition). Comments are preceded by `--`, documentation is enclosed in `/-- ... -/`.

$[\cdot,\cdot]\colon \mathfrak{g} \times \mathfrak{g} \to \mathfrak{g}$ satisfying, for all $X, Y, Z \in \mathfrak{g}$,

$$[X, X] = 0 \tag{1}$$
$$[X, [Y, Z]] = [[X, Y], Z] + [Y, [X, Z]]. \tag{2}$$

The **alternating property** (1) implies skew-symmetry $[X, Y] = -[Y, X]$, and except from characteristic two is equivalent to it. The **Jacobi identity**

$$[X, [Y, Z]] + [Y, [Z, X]] + [Z, [X, Y]] = 0 \tag{3}$$

is a familiar equivalent phrasing of the **Leibnitz rule** (2).

Mathlib has the definition and basic theory of Lie algebras [Nas22]. Conditions (1) and (2) without the vector space structure on $\mathfrak{g}$ and without the bilinearity of the bracket define Mathlib's `LieRing g`, and `LieAlgebra F g` extends this more rudimentary structure. Mathlib's notation for the Lie bracket is `{X, Y}`.

Lie algebra homomorphisms are defined in the usual way, as linear maps which respect the Lie bracket. Mathlib provides the notation `g →₁{F} h` for the type of homomorphisms between two Lie algebras $\mathfrak{g}$ and $\mathfrak{h}$ over $\mathbb{F}$. For comparison, the type of mere $\mathbb{F}$-linear maps is denoted `g →₁[F] h`. A similarity of notation is mostly desirable, but it is worth a warning that this small typographical difference carries substantial content, crucial also in the present report.

Mathlib also has the universal enveloping algebra $\mathcal{U}(\mathfrak{g})$, `UniversalEnvelopingAlgebra F g`. To wit, $\mathcal{U}(\mathfrak{g})$ is an associative algebra over $\mathbb{F}$ together with an inclusion $\iota\colon \mathfrak{g} \hookrightarrow \mathcal{U}(\mathfrak{g})$ which is a Lie algebra homomorphism[2], such that any Lie algebra homomorphism $\phi\colon \mathfrak{g} \to A$ to an associative algebra $A$ over $\mathbb{F}$ lifts uniquely to an algebra homomorphism $\tilde{\phi}\colon \mathcal{U}(\mathfrak{g}) \to A$ which makes the following diagram commute:

$$\begin{array}{ccc} \mathcal{U}(\mathfrak{g}) & & \\ {\scriptstyle \iota}\uparrow & \searrow^{\exists! \tilde{\phi}} & \\ \mathfrak{g} & \xrightarrow{\phi} & A \end{array} \tag{4}$$

A representation of a Lie algebra $\mathfrak{g}$ over $\mathbb{F}$ on a vector space $V$ over $\tilde{\mathbb{F}}$ (with $\tilde{\mathbb{F}}$ possibly extending $\mathbb{F}$) is an $\mathbb{F}$-Lie algebra homomorphism $\rho\colon \mathfrak{g} \to \mathrm{End}_{\tilde{\mathbb{F}}}(V)$. By the universal property of the universal enveloping algebra, this is equivalent to $V$ being a (left) module over the associative algebra $\mathcal{U}(\mathfrak{g})$. In Lean, too, both spellings are possible. We comment on related formalization choices in Sections 11 and 13.

## 3 Central extensions of Lie algebras

The **centre** of a Lie algebra $\mathfrak{e}$ is by definition the set of elements $Z \in \mathfrak{e}$ such that $[Z, E]_\mathfrak{e} = 0$ for all $E \in \mathfrak{e}$.

A Lie algebra $\mathfrak{a}$ is said to be **abelian** if its Lie bracket is identically zero, $[\cdot, \cdot]_\mathfrak{a} = 0$, or equivalently if its centre is the whole of $\mathfrak{a}$. In order to specify an abelian Lie algebra $\mathfrak{a}$, it suffices to give its underlying vector space. The case $\mathfrak{a} = \mathbb{F}$ will be of particular interest for our purposes.[3]

Let $\mathfrak{g}$ be a Lie algebra and let $\mathfrak{a}$ be an abelian Lie algebra.

**Definition 3.1** (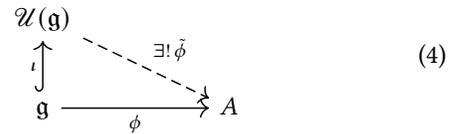). A **central extension** of $\mathfrak{g}$ by $\mathfrak{a}$ is a Lie algebra $\mathfrak{e}$ together with a short exact sequence

$$0 \longrightarrow \mathfrak{a} \xhookrightarrow{\iota} \mathfrak{e} \xtwoheadrightarrow{\pi} \mathfrak{g} \longrightarrow 0 \tag{5}$$

of Lie algebras such that $\iota(\mathfrak{a})$ is contained in the centre of $\mathfrak{e}$.

The mapping $\iota$ provides an (injective) inclusion of a copy $\iota(\mathfrak{a})$ of $\mathfrak{a}$ in (the centre of) the central extension $\mathfrak{e}$, whereas the mapping $\pi$ provides a (surjective) projection from $\mathfrak{e}$ to $\mathfrak{g}$. Exactness at the middle requires $\mathrm{Ker}(\pi) = \iota(\mathfrak{a})$, so $\pi$ factors to give an isomorphism $\mathfrak{e}/\iota(\mathfrak{a}) \cong \mathfrak{g}$.

---

[2] When needed, we equip associative algebras with their natural Lie bracket given by $[a, b] = ab - ba$. This is also a readily available `instance` in Mathlib, as discussed in [Nas22, Sec. 5].

[3] The main reason is that projective representations of a Lie algebra $\mathfrak{g}$ correspond to (ordinary) representations of central extensions of $\mathfrak{g}$ by $\mathbb{F}$, and projective representations are the way symmetries act in quantum physics.

As for when two central extensions should be thought of as "the same", the natural notion is:

**Definition 3.2.** Two central extensions $0 \longrightarrow \mathfrak{a} \xrightarrow{\iota_j} \mathfrak{e}_j \xrightarrow{\pi_j} \mathfrak{g} \longrightarrow 0$, $j = 1, 2$, are said to be **equivalent** if there exists a Lie algebra isomorphism $\phi \colon \mathfrak{e}_1 \to \mathfrak{e}_2$ such that the diagram

$$\begin{array}{c} \mathfrak{e}_1 \\ {}^{\iota_1}\nearrow \quad \downarrow \phi \cong \quad \searrow^{\pi_1} \\ \mathfrak{a} \qquad\qquad \mathfrak{g} \\ {}_{\iota_2}\searrow \quad \uparrow \quad \nearrow_{\pi_2} \\ \mathfrak{e}_2 \end{array} \qquad (6)$$

commutes.

By (5), a central extension $\mathfrak{e}$ of $\mathfrak{g}$ by $\mathfrak{a}$ is, as a vector space, isomorphic to the direct sum $\mathfrak{g} \oplus \mathfrak{a}$. If we fix such a vector space isomorphism, we may consider elements of $\mathfrak{e}$ as pairs $(X, A)$, where $X \in \mathfrak{g}$ and $A \in \mathfrak{a}$, and the maps in (5) take the forms $\iota(A) = (0, A)$ and $\pi((X, A)) = X$ for $X \in \mathfrak{g}$, $A \in \mathfrak{a}$. For $\pi$ in (5) to be a Lie algebra homomorphism, Lie brackets in $\mathfrak{e}$ must then take the form $[(X, A), (Y, B)]_\mathfrak{e} = ([X, Y]_\mathfrak{g}, C)$. Moreover, for $(0, A)$ and $(0, B)$ to be central in $\mathfrak{e}$, these brackets cannot depend on $A$ and $B$, so the second component $C$ is a function of $X$ and $Y$ only, $C = \omega(X, Y)$ for some $\omega \colon \mathfrak{g} \times \mathfrak{g} \to \mathfrak{a}$. For the bracket in $\mathfrak{e}$ to be bilinear and antisymmetric, $\omega$ must be as well. For the Leibniz rule (2) to hold in $\mathfrak{e}$, $\omega$ must satisfy

$$\omega(X, [Y, Z]) = \omega([X, Y], Z) + \omega(Y, [X, Z]) \qquad (7)$$

for all $X, Y, Z \in \mathfrak{g}$.

Conversely, given a bilinear antisymmetric $\omega \colon \mathfrak{g} \times \mathfrak{g} \to \mathfrak{a}$ satisfying (7), we may construct a central extension $\mathfrak{e}_\omega$ by equipping the vector space $\mathfrak{g} \oplus \mathfrak{a}$ with the bracket

$$[(X, A), (Y, B)]_{\mathfrak{e}_\omega} = ([X, Y]_\mathfrak{g}, \omega(X, Y)). \qquad (8)$$

Note, however, that given any linear map $\beta \colon \mathfrak{g} \to \mathfrak{a}$, the map $\phi \colon \mathfrak{g} \oplus \mathfrak{a} \to \mathfrak{g} \oplus \mathfrak{a}$ defined by $\phi((X, A)) = (X, A - \beta(X))$ is a linear isomorphism which makes (6) commute. Since

$$\begin{aligned} \big[\phi\big((X, A)\big), \phi\big((Y, B)\big)\big] &= [(X, A - \beta(X)), (Y, B - \beta(Y))] \\ &= ([X, Y]_\mathfrak{g}, \omega(X, Y)) \\ &= \phi\big(([X, Y]_\mathfrak{g}, \ \omega(X, Y) + \beta([X, Y]_\mathfrak{g}))\big), \end{aligned}$$

this is a Lie algebra isomorphism and an equivalence of central extensions defined by $\omega$ and by

$$(X, Y) \mapsto \omega(X, Y) + \beta([X, Y]_\mathfrak{g}). \qquad (9)$$

All equivalences of central extensions with $\iota \colon A \mapsto (0, A)$ and $\pi \colon (X, A) \mapsto X$ arise in this way.

## 4 Lie algebra cohomology in degree two

The considerations in Section 3 give a systematic way to classify central extensions of $\mathfrak{g}$ by $\mathfrak{a}$, in cohomological terms. Lie algebra cohomology generally was introduced in 1948 by Chevalley and Eilenberg [CE48]; the following are the special case definitions in degree two, directly suggested by the study of central extensions.

Considerations of constructing central extensions suggest:

**Definition 4.1 (↗).** A 2-**cocycle** of $\mathfrak{g}$ with coefficients in $\mathfrak{a}$ is a bilinear antisymmetric map $\omega \colon \mathfrak{g} \times \mathfrak{g} \to \mathfrak{a}$ satisfying (7). Denote the vector space of such 2-cocycles by $Z^2(\mathfrak{g}; \mathfrak{a})$.

Considerations of equivalence of central extensions, particularly (9), further suggest:

**Definition 4.2 (↗).** A 1-**cochain** of $\mathfrak{g}$ with coefficients in $\mathfrak{a}$ is a linear map $\beta \colon \mathfrak{g} \to \mathfrak{a}$. Denote the vector space of such 1-cochains by $C^1(\mathfrak{g}; \mathfrak{a})$. For $\beta \in C^1(\mathfrak{g}; \mathfrak{a})$, define $\partial \beta \colon \mathfrak{g} \times \mathfrak{g} \to \mathfrak{a}$ by

$$\partial \beta(X, Y) = \beta([X, Y]_\mathfrak{g}), \qquad (10)$$

which by bilinearity, antisymmetry, and Leibniz rule (2) is a 2-cocycle, $\partial \beta \in Z^2(\mathfrak{g}; \mathfrak{a})$. Denote the range of the linear map $\partial \colon C^1(\mathfrak{g}; \mathfrak{a}) \to Z^2(\mathfrak{g}; \mathfrak{a})$ by $B^2(\mathfrak{g}; \mathfrak{a})$, and call its elements 2-**coboundaries**.

For systematically classifying central extensions, we need:

**Definition 4.3 (↗).** The 2-**cohomology** of $\mathfrak{g}$ with coefficients in $\mathfrak{a}$ is the quotient space

$$H^2(\mathfrak{g}; \mathfrak{a}) \ = \ Z^2(\mathfrak{g}; \mathfrak{a})/B^2(\mathfrak{g}; \mathfrak{a}). \qquad (11)$$

The classification result which follows from the considerations sketched in Section 3 is:

**Theorem 4.4.** *Let $\mathfrak{g}$ be a Lie algebra and $\mathfrak{a}$ an abelian Lie algebra. Any 2-cocycle $\omega \in Z^2(\mathfrak{g}; \mathfrak{a})$ defines a central extension $\mathfrak{e}_\omega$ of $\mathfrak{g}$ by $\mathfrak{a}$, and two cocycles $\omega_1, \omega_2 \in Z^2(\mathfrak{g}; \mathfrak{a})$ define equivalent central extensions if and only if $\omega_2 - \omega_1 \in B^2(\mathfrak{g}; \mathfrak{a})$. Moreover, any central extension is equivalent to a central extension defined by a 2-cocycle (unique modulo $B^2(\mathfrak{g}; \mathfrak{a})$).*

The zero cocycle $0 \in Z^2(\mathfrak{g}; \mathfrak{a})$ defines the trivial central extension: the direct sum Lie algebra $\mathfrak{g} \oplus \mathfrak{a}$ (the vector space direct sum equipped with the direct sum Lie bracket). Cocycles $\omega \in Z^2(\mathfrak{g}; \mathfrak{a})$ whose cohomology class $[\omega] \in H^2(\mathfrak{g}; \mathfrak{a})$ is zero define central extensions equivalent to this trivial one. Nontrivial central extensions are obtained from cocycles $\omega$ which are cohomologically nontrivial, $[\omega] \neq 0$.

A prominent class of Lie algebras is that of the finite-dimensional complex semisimple Lie algebras, see [Bou02] and [Nas22, Sec. 10]. Already the seminal paper [CE48] noted that the second cohomology of such a Lie algebra vanishes, implying the absence of nontrivial central extensions. Interesting examples, which still share many properties with semisimple Lie algebras, arise in infinite dimensions. From Section 6 on, we focus on one such case.

First, however, we pause to comment on the formalization of the above generalities.

## 5 Formalizing central extensions

A general formalization principle distilled from experience is to systematically separate two aspects of what in informal mathematics often get conflated as "a definition": a *characteristic predicate* encodes the defining property part of a definition, and a *construction* provides a concrete object with that property, which is then meant to be used subsequently as "the defined object". Sometimes it is useful to provide multiple alternative constructions (whose isomorphism typically follows by the characteristic predicate), but usually the precise construction does not matter, as we simply want an object with the desired properties.

The discussion on central extensions in Section 3 readily suggests such a division: Definition 3.1 is a defining property, whereas formula (8) yields a concrete construction.[4]

The Lean setup to start discussing central extensions is:

```
variable {F : Type*} [CommRing F]
variable {g : Type*} [LieRing g] [LieAlgebra F g]
variable {a : Type*} [LieRing a] [LieAlgebra F a]
variable {e : Type*} [LieRing e] [LieAlgebra F e]
```

In particular, for now it is enough to assume that $\mathbb{F}$ is a commutative ring instead of a field; formalization benefits from minimizing hypotheses at all stages. Definition 3.1 is then formally implemented as a predicate on the pair $\iota\colon \mathfrak{a} \to \mathfrak{e}$, $\pi\colon \mathfrak{e} \to \mathfrak{g}$ of Lie algebra homomorphisms:

```
structure LieAlgebra.IsExtension
    (i : a →₁{F} e) (p : e →₁{F} g) : Prop where
  ker_eq_bot : i.ker = ⊥    -- injectivity of `i`
  range_eq_top : p.range = ⊤ -- surjectivity of `p`
  exact : i.range = p.ker    -- exactness

structure LieAlgebra.IsCentralExtension
    (i : a →₁{F} e) (p : e →₁{F} g)
    extends IsExtension i p where
  central : ∀ (A : a), ∀ (E : e), {i A, E} = 0
```

As observed in Section 3, this defining property should let us identify $\mathfrak{e}$ with $\mathfrak{g} \oplus \mathfrak{a}$ as vector spaces. The identification is not canonical, it depends on a choice of a "section" $\sigma\colon \mathfrak{g} \to \mathfrak{e}$ of $\pi$ (a linear map satisfying $\pi \circ \sigma = \mathrm{id}_{\mathfrak{g}}$). With a section specified, one can in particular construct a basis of $\mathfrak{e}$ from bases of $\mathfrak{g}$ and $\mathfrak{a}$; this will be later used for the Virasoro algebra and the Heisenberg algebra:

```
noncomputable def LieAlgebra.IsExtension.basis
    {i : a →₁{F} e} {p : e →₁{F} g}
    (ex : LieAlgebra.IsExtension i p)
    (σ : g →₁{F} e) (hσ : p.toLinearMap ∘₁ σ = 1)
    {ιA ιG : Type*}
    (basA : Basis ιA F a) (basG : Basis ιG F g) :
    Basis (ιA ⊕ ιG) F e := ...
```

For the construction part, the setup is:

```
variable (F : Type*) [CommRing F]
variable (g : Type*) [LieRing g] [LieAlgebra F g]
variable (a : Type*) [AddCommGroup a] [Module F a]
```

Definitions 4.1 and 4.2 are phrased in Lean as the types:

```
structure LieOneCochain where
  toLinearMap : g →₁[F] a

structure LieTwoCocycle where
  toBilin : g →₁[F] g →₁[F] a
  self' : ∀ X, toBilin X X = 0
  leibniz' : ∀ X Y Z,
    toBilin X {Y, Z}
      = toBilin {X, Y} Z + toBilin Y {X, Z}
```

We note here that bilinear maps are written in Lean as linear maps to a space of linear maps, so `toBilin` in `LieTwoCocycle` above indeed is the data of a bilinear map $\mathfrak{g} \times \mathfrak{g} \to \mathfrak{a}$.

It is routine to equip the types `LieTwoCocycle F g a` and `LieOneCochain F g a` with the vector space structure, i.e., provide instances of `AddCommGroup` and `Module F` for them. A more interesting Mathlib design is to provide `LinearMapClass` instances which allow interpreting terms of these types as (bi)linear maps, and which directly give access to a well-developed API of (bi)linear maps. Such designs are one way in which Mathlib avoids code duplication and ensures unity of the library.

The linear map $\partial\colon C^1(\mathfrak{g}, \mathfrak{a}) \to Z^2(\mathfrak{g}, \mathfrak{a})$ is a term of type `LieOneCochain F g a →₁[F] LieTwoCocycle F g a`, which we named `LieOneCochain_bdryHom`. We also provide `LieOneCochain.bdry` to enable convenient dot notation such as `β.bdry`.[5] Coboundaries and cohomology are then defined, corresponding to Definitions 4.2 and 4.3, as:

```
abbrev LieTwoCoboundary :=
  LinearMap.range (LieOneCochain_bdryHom F g a)

def LieTwoCohomology :=
  LieTwoCocycle F g a / LieTwoCoboundary F g a
```

For constructing central extensions from 2-cocycles, we use another formalization idiom: *type synonyms*. We want to use the underlying direct sum vector space[6] $\mathfrak{g} \oplus \mathfrak{a}$, but we want to avoid using the direct sum Lie algebra structure. Therefore we create an opaque type

```
def LieTwoCocycle.CentralExtension
    (ω : LieTwoCocycle F g a) :=
  g × a
```

---

[4] Here the need to separate the predicate from the construction is rather evident, since Definition 3.1 does not by itself completely characterize a central extension; the additional data of the cohomology class is required to get uniqueness up to isomorphism. A perhaps more representative example of the same formalization principle appears in Section 13.

[5] The reason to have both is that this dot notation does not allow bundling the linear dependence on `β : LieOneCochain F g a`.

[6] The most straightforward spelling of the direct sum of just two vector spaces is the type `g × a`, which by standard typeclass instances comes equipped with the direct sum vector space structure.

and equip it with the appropriate type class instances, including a `LieRing` and `LieAlgebra` instance with the Lie bracket depending on the cocycle $\omega$ as in (8): ↗

```
instance : LieRing ω.CentralExtension where
  bracket Z W := ⟨{Z.fst, W.fst}, ω Z.fst W.fst⟩
  ... -- omit properties of bracket for brevity
```

Naturally, then, we should show that the constructed central extension satisfies the defining property of central extensions. This requires introducing the Lie algebra homomorphisms corresponding to $\iota\colon \mathfrak{a} \to \mathfrak{g} \oplus \mathfrak{a}$ and $\pi\colon \mathfrak{g} \oplus \mathfrak{a} \to \mathfrak{g}$:

```
variable (ω : LieTwoCocycle F g a)

def LieTwoCocycle.CentralExtension.emb :
    a →₁{F} ω.CentralExtension where
  toFun := fun A ↦ ⟨0, A⟩
  ... -- omit homomorphism properties for brevity

def LieTwoCocycle.CentralExtension.proj :
    ω.CentralExtension →₁{F} g where
  toFun := fun ⟨X, _⟩ ↦ X
  ... -- omit homomorphism properties for brevity
```

For example the formal statement that a central extension constructed from a cocycle $\omega$ satisfies the characteristic predicate of extensions is ↗

```
instance LieTwoCocycle.CentralExtension.isExtension :
    LieAlgebra.IsExtension (emb ω) (proj ω) := ...
```

for which we just supply the trivial proofs of `(emb ω).ker = ⊥`, `(proj ω).range = ⊤`, `(emb ω).range = (proj ω).ker`.

For the purpose of constructing obvious bases for extensions, we also introduce the obvious section $\sigma\colon \mathfrak{g} \to \mathfrak{g} \oplus \mathfrak{a}$, $X \mapsto (X, 0)$ under the name `stdSection`.

One mild inconvenience arising from dependent type theory may be worth mentioning. The type `ω.CentralExtension` depends on the cocycle `ω : LieTwoCocycle F g a` (that is its whole point!). Recall from Section 3 and (9) the essential property that two cocycles differing by a coboundary define equivalent central extensions. So for $\omega \in Z^2(\mathfrak{g}; \mathfrak{a})$ and $\beta \in C^1(\mathfrak{g}; \mathfrak{a})$, we must exhibit a Lie algebra isomorphism $\mathfrak{e}_\omega \cong \mathfrak{e}_{\omega+\partial\beta}$. First implement, straightforwardly, a Lie algebra homomorphism $\mathfrak{e}_\omega \to \mathfrak{e}_{\omega+\partial\beta}$: ↗

```
def LieOneCochain.bdryHom (β : LieOneCochain F g a)
    (ω : LieTwoCocycle F g a) :
    (ω.CentralExtension)
      →₁{F} (ω + β.bdry).CentralExtension := ...
```

Mathematically, the inverse of this homomorphism is of precisely the same form, $\mathfrak{e}_{\omega+\partial\beta} \to \mathfrak{e}_\omega$ — just start from the cocycle $\omega + \partial\beta \in H^2(\mathfrak{g}; \mathfrak{a})$ and modify by the opposite coboundary $-\partial\beta \in B^2(\mathfrak{g}; \mathfrak{a})$. But in dependent type theory, the composition lands in $\mathfrak{e}_{\omega+\partial\beta-\partial\beta}$, which is not definitionally equal to $\mathfrak{e}_\omega$. In this case the workaround seems to be to first prove that if we have a propositional equality $\omega_1 = \omega_2$ (an equality that is provable but not necessarily definitional), then we have a Lie algebra isomorphism $\mathfrak{e}_{\omega_1} \cong \mathfrak{e}_{\omega_2}$: ↗

```
def LieTwoCocycle.CentralExtension.congr
    {ω₁ ω₂ : LieTwoCocycle F g a} (h : ω₁ = ω₂) :
    ω₁.CentralExtension
      ≃l{F} ω₂.CentralExtension := ...
```

Half a dozen small lemmas need to be added about this `LieTwoCocycle.CentralExtension.congr`, and then one has the desired isomorphism: ↗

```
noncomputable def equiv_of_lieTwoCoboundary
    {ω ω' : LieTwoCocycle F g a}
    (h : ω' - ω ∈ LieTwoCoboundary F g a) :
    (ω.CentralExtension) ≃l{F} (ω'.CentralExtension)
```

Thus the problem is not serious, but the dependent type issue at its core is similar to what held back the development of practically usable homological algebra in Lean for some time; see [Liv23, Sec. 2.2.1] for a discussion.

## 6 The Witt algebra

The Witt algebra is essentially the Lie algebra of smooth vector fields on the circle, or the Lie algebra of holomorphic vector fields in an annular region; see Figure 2 for illustrations. Since flows of holomorphic vector fields are conformal transformations, the Witt algebra can be regarded as the Lie algebra of infinitesimal conformal transformations in two dimensions. This interpretation explains its relevance to conformal field theory.

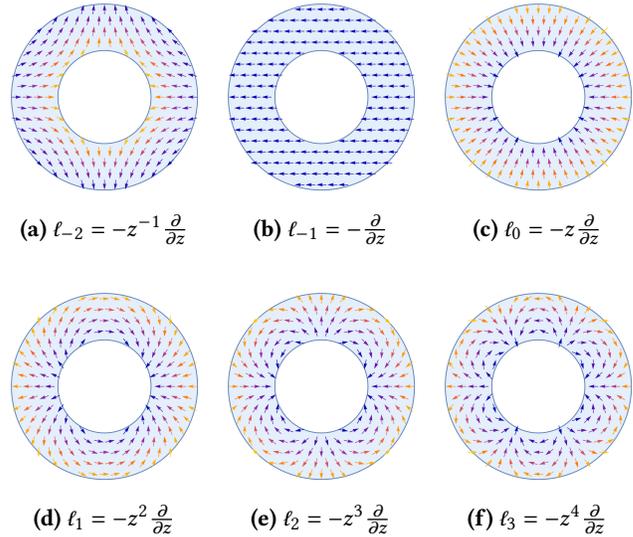

**(a)** $\ell_{-2} = -z^{-1}\frac{\partial}{\partial z}$  **(b)** $\ell_{-1} = -\frac{\partial}{\partial z}$  **(c)** $\ell_0 = -z\frac{\partial}{\partial z}$

**(d)** $\ell_1 = -z^2\frac{\partial}{\partial z}$  **(e)** $\ell_2 = -z^3\frac{\partial}{\partial z}$  **(f)** $\ell_3 = -z^4\frac{\partial}{\partial z}$

**Figure 2.** The Witt algebra $\mathfrak{witt}$ can be viewed as the Lie algebra of (polynomial) vector fields on the circle, or as the Lie algebra of (Laurent polynomial) holomorphic vector fields in an annular region.

The definition, however, is most straightforwardly given in terms of a concrete basis $(\ell_n)_{n \in \mathbb{Z}}$. In the interpretation as holomorphic vector fields in $\mathbb{C} \setminus \{0\}$, the basis element $\ell_n$ corresponds to the vector field $-z^{n+1}\frac{\partial}{\partial z}$.

**Definition 6.1 (↗).** Fix a ground field $\mathbb{F}$. The **Witt algebra** $\mathfrak{witt}$ is the Lie algebra over $\mathbb{F}$ with a basis $(\ell_n)_{n \in \mathbb{Z}}$ and the

Lie bracket given in this basis as

$$[\ell_m, \ell_n] = (m-n)\ell_{m+n}. \tag{12}$$

According to the general principle (Theorem 4.4), central extensions of $\mathfrak{witt}$ are classified by the second Lie algebra cohomology group $H^2(\mathfrak{witt}; \mathbb{F})$. This and the entire cohomology ring $H^*(\mathfrak{witt}; \mathbb{F})$, was first calculated by Gelfand and Fuks in 1969 [GF69].[7]

**Theorem 6.2** ( [GF69]). *Suppose that $\mathbb{F}$ is of characteristic 0. Then:*

(a) *The formula*
$$\omega_{\mathfrak{vir}}(\ell_m, \ell_n) = \frac{m^3 - m}{12} \delta_{m,-n} \tag{13}$$
*defines a 2-cocycle of $\mathfrak{witt}$ with coefficients in $\mathbb{F}$.*
(b) *The cohomology class $[\omega_{\mathfrak{vir}}] \in H^2(\mathfrak{witt}; \mathbb{F})$ is nonzero.*
(c) *The cohomology group $H^2(\mathfrak{witt}; \mathbb{F})$ is one-dimensional and generated by $[\omega_{\mathfrak{vir}}]$*
$$H^2(\mathfrak{witt}; \mathbb{F}) = \mathrm{span}_{\mathbb{F}}\{[\omega_{\mathfrak{vir}}]\}.$$

The proofs of (a), (b), (c) are calculations, included in Appendix A.1 and in the online blueprint.

## 7 Formalizing Witt algebra cohomology

The very definition of a Lie algebra involves a trilinear identity: the Jacobi identity (3), or equivalently the Leibniz rule (2). Similarly the definition of a Lie algebra 2-cocycle (7) involves a trilinear identity. When constructing concrete Lie algebras or Lie algebra cocycles, e.g. Definition 6.1 and (13), one typically verifies the equations in a basis and appeals to multilinearity. However obvious mathematically, the first step formally is to argue that the relation to be established is multilinear. We introduced `cyclicTripleSumHom` which constructs a trilinear map

$$(X, Y, Z) \mapsto \mu(X, \nu(Y,Z)) + \mu(Y, \nu(Z,X)) + \mu(Z, \nu(X,Y))$$

out of bilinear maps $\nu \colon V \times V \to V$ and $\mu \colon V \times V \to W$. Proving trilinearity amounts to proving additivity and scalar multiplication property in each of the three inputs; cyclicity allows to reduce to one variable. For Jacobi identity one uses $V = W = \mathfrak{g}$ and for 2-cocycle condition $V = \mathfrak{g}$ and $W = \mathbb{F}$; this reduced code duplication in the formalization of Definition 6.1 and (13), in particular. We expect these patterns to be occasionally useful also elsewhere in the formalization of Lie theory.

In the formalization of Theorem 6.2, the standing assumptions about the field are now:

```
variable (F : Type*) [Field F] [CharZero F]
```

---
[7]The article [GF69] (in Russian) slightly predates [Vir70]. We use the term *Virasoro cocycle* for (13) (instead of, e.g., *Feigin-Fuks cocycle*) merely because of its use is Definition 8.1. For an account of the history and origins of the cocycle, see [Rog10].

Part (a) of Theorem 6.2 is in fact phrased as a definition `WittAlgebra.virasoroCocycle F` of type described in Section 5, `LieTwoCocycle F (WittAlgebra F) F`. A statement demonstrating that the definition is as intended in (13)
```
lemma virasoroCocycle_apply_lgen_lgen (n m : ℤ) :
    virasoroCocycle F (lgen F n) (lgen F m)
      = if n + m = 0 then (n^3 - n : F) / 12 else 0
```
is recorded, which uses an explicit basis `WittAlgebra.lgen F` of the Witt algebra. Parts (b) and (c) are stated in Lean as:
```
theorem cohomologyClass_virasoroCocycle_ne_zero :
    (virasoroCocycle k).cohomologyClass ≠ 0
theorem WittAlgebra.rank_lieTwoCohomology_eq_one :
    Module.rank F
      (LieTwoCohomology F (WittAlgebra F) F) = 1
```
The proofs are largely calculations (see Appendix A.1), and as such they differ from the majority of mathematics that is formalized, e.g., in Mathlib. We discuss the experience of writing and maintaining calculational proofs in Section 11, which deals with more substantial calculations.

## 8 The Virasoro algebra

After Definition 6.1 and Theorem 6.2 and the principle encoded by Theorem 4.4, the definition of the Virasoro algebra is straightforward.[8]

**Definition 8.1** (). The **Virasoro algebra** $\mathfrak{vir}$ over a field $\mathbb{F}$ of characteristic 0 is the central extension of the Witt algebra $\mathfrak{witt}$ associated to the 2-cocycle $\omega_{\mathfrak{vir}}$ in (13).

The definition of $\mathfrak{vir}$ is equally straighforwad in Lean:
```
def VirasoroAlgebra
    (F : Type*) [Field F] [CharZero F] :=
  LieTwoCocycle.CentralExtension
    (WittAlgebra.virasoroCocycle F)
```
The Lie algebra structure on the type `VirasoroAlgebra F` is obtained immediately through the general constructions of Section 5, `LieTwoCocycle.CentralExtension.instLieAlgebra`.

From Definition 8.1 (along with Definitions 6.1 and 3.1) it is clear that $\mathfrak{vir}$ has a basis consisting of $\mathcal{L}_n$, $n \in \mathbb{Z}$, and a central element $C$, such that the Lie brackets of these basis elements are

$$[\mathcal{L}_m, \mathcal{L}_n] = (m-n)\mathcal{L}_{m+n} + \frac{m^3 - m}{12}\delta_{m,-n} C \tag{14}$$
$$[C, \mathcal{L}_n] = 0.$$

The basis element $C$ is central, and typically we study representations of $\mathfrak{vir}$, where it acts as a scalar multiple of the identity operator.[9] The value $c \in \mathbb{F}$ of that scalar is then called the **central charge** of the representation.

To get this basis in Lean, we use `stdSection` and `LieAlgebra.IsExtension.basis` from Section 5. We also provide the notation `lgen F n` and `cgen F` for $\mathcal{L}_n$ and $C$ in the `VirasoroAlgebra` namespace. The statements asserting that the basis is as in (14) are then:

```
lemma lgen_bracket (n m : ℤ) :
    ⁅lgen F n, lgen F m⁆
      = (n - m : F) • lgen F (n + m)
        + if n + m = 0
            then ((n^3 - n : F)/12) • cgen F else 0

lemma cgen_bracket (X : VirasoroAlgebra F) :
    ⁅cgen F, X⁆ = 0
```

---
[8]The Lie algebra $\mathfrak{vir}$ is nowadays referred to as the Virasoro algebra, owing to a series of impactful works that started from [Vir70].
[9]In an irreducible representation, for example, this is necessarily the case.

## 9 The Heisenberg algebra

A parallel but simpler story allows to define the Heisenberg algebra. This is the simplest example of current algebras, relevant to the conformal field theory of the free boson field, and a starting point for several more intricate constructions in vertex operator algebras and conformal field theory.

For this purpose, start from the abelian Lie algebra with basis $(\mathcal{I}_k)_{k \in \mathbb{Z}}$, and form a central extension using its (nontrivial) 2-cocycle given by the formula $\omega_{\mathfrak{hei}}(\mathcal{I}_k, \mathcal{I}_l) = k\, \delta_{k,-l}$. This gives rise to the **Heisenberg algebra** $\mathfrak{hei}$, a Lie algebra with a basis consisting of $\mathcal{J}_k$, $k \in \mathbb{Z}$, and a central element $\mathcal{K}$, with brackets

$$[\mathcal{J}_k, \mathcal{J}_l] = k\, \delta_{k,-l}\, \mathcal{K} \qquad (15)$$
$$[\mathcal{K}, \mathcal{J}_k] = 0.$$

The basis element $\mathcal{K}$ is central, and typically we study representations of $\mathfrak{hei}$, where it acts as a scalar multiple of the identity operator. If that scalar is nonzero, in view of (15), we may normalize the basis so that the value of the scalar is 1, i.e., $\mathcal{K}$ acts simply as the identity.

## 10 The basic Sugawara construction

There are a number of variants of Sugawara constructions. The simplest is a (physically motivated) way to equip a representation of the Heisenberg algebra subject to some hypotheses with the structure of a representation of the Virasoro algebra with central charge $c = 1$. We call that particular variant the **basic bosonic Sugawara construction**, and we outline it below. Typically in the literature, such constructions are presented for very particular representations (e.g., highest weight representations of affine Lie algebras or modules over vertex algebras) [GO85, Kac98, LL04], but upon inspection of the proofs it is clear that modest and easily stated hypotheses suffice. The generality gained by working under minimal hypotheses is genuinely important for example in applications to local fields of probabilistic lattice models [HKV22, AC24, ABK24], where very little information is a priori available about the structure of the representation one has.

Fix a field $\mathbb{F}$ of characteristic 0. Let $V$ be a vector space, and let $(\mathsf{J}_k)_{k \in \mathbb{Z}}$ be a collection of linear operators $\mathsf{J}_k \colon V \to V$ about which we make the following two assumptions:

(H) Assume the Heisenberg algebra relations (15) for these operators in the following form:[10]

$$\mathsf{J}_k \circ \mathsf{J}_l - \mathsf{J}_l \circ \mathsf{J}_k = k\, \delta_{k,-l}\, \mathrm{id}_V \qquad \text{for } k, l \in \mathbb{Z}. \qquad (16)$$

(T) Assume local truncation: for any $v \in V$ there exists a $K_v \in \mathbb{Z}$ such that

$$\mathsf{J}_l v = 0 \qquad \text{for all } l \geq K_v. \qquad (17)$$

We first define normal ordered quadratic expressions in the operators $(\mathsf{J}_k)_{k \in \mathbb{Z}}$ by

$${:}\mathsf{J}_k \mathsf{J}_l{:} \; = \begin{cases} \mathsf{J}_k \circ \mathsf{J}_l & \text{if } k \leq l, \\ \mathsf{J}_l \circ \mathsf{J}_k & \text{if } k > l. \end{cases} \qquad (18)$$

Using these, we define operators $\mathsf{L}_n \colon V \to V$, $n \in \mathbb{Z}$, by

$$\mathsf{L}_n v = \frac{1}{2} \sum_{k \in \mathbb{Z}} {:}\mathsf{J}_{n-k} \mathsf{J}_k{:} v \qquad \text{for } v \in V, \qquad (19)$$

noting that the formally infinite sum contains only finitely many nonzero terms by virtue of the local truncation hypothesis (T) and the definition (18). The number of terms needed, however, depends on the vector $v \in V$ on which the operator $\mathsf{L}_n$ is applied, so the operators $\mathsf{L}_n$ themselves are not given by finite sums of (18).[11]

The basic bosonic Sugawara construction is the following.

**Theorem 10.1** (). *Under the hypotheses (H) and (T) above, the linear operators* (19) *on $V$ satisfy the commutation relations*

$$\mathsf{L}_m \circ \mathsf{L}_n - \mathsf{L}_n \circ \mathsf{L}_m = (m - n)\, \mathsf{L}_{m+n} + \frac{m^3 - m}{12}\, \delta_{m,-n}\, \mathrm{id}_V.$$

**Corollary 10.2** (). *Under the hypotheses (H) and (T), a representation of the Virasoro algebra on $V$, with central charge $c = 1$, is defined by the linear extension of the assignments $\mathcal{L}_n \mapsto \mathsf{L}_n$, $C \mapsto \mathrm{id}_V$.*

The proof of Theorem 10.1 is a calculation. The details of it are included in the online blueprint and in Appendix A.2.

## 11 Formalizing Sugawara constructions

The main design decision worth noting is the handling of apparently infinite sums such as (19). The indexing type `ℤ` is infinite, but any time an operator $\mathsf{L}_n$ is acting on a vector $v \in V$, only finitely many of the terms are nonzero by (T).

A standard way to handle infinite sums is to have a notion of their convergence; Mathlib's topological sum notion is `tsum` (`∑'`). We could indeed equip $V$ with discrete topology, and prove that convergence of an infinite sum, expressed as the predicate `Summable`, is equivalent with only finitely

---

[10]This corresponds to (15) if the central element $\mathcal{K}$ is normalized to act as the identity.

[11]The literature (including Wikipedia at the time of writing) contains statements such as *"The Sugawara construction is an embedding of the Virasoro algebra into the universal enveloping algebra of the affine Lie algebra"*, which are false for the same reason, although they may be helpful to convey a rough idea of what the Sugawara construction does.

Virasoro algebra and Sugawara constructions formally in Lean

many of its terms being nonzero. In arguments about the sums, one would need to establish and carry around some such summability hypotheses. The main reason to discard this design, however, is the artificial and forceful placing of the discrete topology on $V$. Mathematically we often do not want to put any topology on $V$, or worse, we might want to equip $V$ with an altogether different topology!

Another alternative is to fully embrace the finiteness of the sums, and write them using Mathlib's `Finset.sum` ( $\Sigma$ ) by explicitly specifying which finite subset of indices is used at each step. Since the condition (T) is not quantitative, to get the finite sets requires invoking some version of the axiom of choice for every vector that appears, even in intermediate steps of a calculation. This hands-on bookkeeping approach becomes very tedious, while being totally beside the point of proving the interesting commutation relations!

Mathlib offers a convenient design suitable for the situation: `finsum` ( $\Sigma^f$ ). This is a notion of sums which Lean allows the user to write down without topology and without an a priori finite support: it is just defined to yield *junk values* (zero) in case a user were to unwisely consider a sum with infinitely many nonzero terms. In arguments manipulating the sums, sometimes finite support hypotheses are needed, just like the summability hypotheses are needed for topological sums. But as opposed to the constructivist feeling with `Finset.sum`, proving finiteness when such hypotheses are needed aligns well with (careful) mathematical reasoning. Most importantly, `finsum` eliminates the notational and formalization burden of explicitly providing a set and a proof of its finiteness each time a sum is mentioned.

The setup of operators $(\mathsf{J}_k)_{k \in \mathbb{Z}}$ on a vector space $V$ is:

```
variable {F : Type*} [Field F]
variable {V : Type*} [AddCommGroup V] [Module F V]
variable (heiOper : ℤ → (V →₁[F] V))
```

The hypothesis (H) is written as[12]

```
variable (heiComm : ∀ k l,
         (heiOper k).commutator (heiOper l)
           = if k + l = 0 then (k : F) • 1 else 0)
```

The hypothesis (T) is idiomatically phrased using the filter `atTop` on the set $\mathbb{Z}$ indexing the collection `heiOper`:

```
variable (heiTrunc : ∀ v,
         atTop.Eventually fun l ↦ heiOper l v = 0)
```

Constructing the operators (19) as linear self-maps of $V$, i.e., as terms of type `V →₁[F] V`, already involves ensuring that the formally infinite sums can be defined without breaking properties such as additivity. Therefore the definition `sugawaraGen heiTrunc n` corresponding to $\mathsf{L}_n$ needs the argument `heiTrunc` encoding the truncation property (T). The statement of the final result, Theorem 10.1, is:

```
lemma commutator_sugawaraGen [CharZero F] (n m : ℤ) :
    (sugawaraGen heiTrunc n).commutator
      (sugawaraGen heiTrunc m)
    = (n - m) • (sugawaraGen heiTrunc (n + m))
      + if n + m = 0 then
          ((n ^ 3 - n : F) / (12 : F)) • 1
        else 0
```

Again, the formal proof of Theorem 10.1 follows a fairly straightforward translation of the calculation needed in the informal mathematical proof (see Appendix A.2). We note that a stronger integration of symbolic calculation to the Lean theorem prover has the potential to simplify proofs of such calculational flavor. The existing tactics `ring`, `simp`, `norm_num`, `congr` were useful in automating and structuring (steps of) calculation. The release of the powerful `grind` tactic in the middle of the present project resulted in a few additional simplifications. If the Gröbner basis based algebraic simplification capabilities of `grind` were suitably extended with on Gröbner–Shirshov bases, one could expect very useful automation of calculations in the universal enveloping algebras of Lie algebras.

When updating versions of Lean and Mathlib, calculational proofs tended to be more prone to breaking than other proofs. One explanation is that changes, e.g., to `simp`-tagged lemmas and refactors in the library cause the behavior of the automation to change. But we suspect that the brittleness was partly a fault of the way our proofs were written.

Specifically, we expect that calculational proofs arranged to simple `calc` steps should be significantly easier to fix during updates and refactors, but most of our proofs were written instead step by step in the tactic mode.[13] Given that `calc` proofs are also more readable, they seem like an overall better organization of calculational proofs. Their only downside is the process of constructing a formal proof: making progress by manipulating the tactic state feels easier. The community seems aware of this tension, and Lean already provides tools to facilitate writing `calc` proofs. With the arrival of Lean 4, widgets were introduced to manipulate `calc` state in a manner very natural to a mathematician, but in practice these are still not easy to adopt as a part of one's workflow. The recent `calcify` tactic can turn a proof written in tactic mode (or term mode) to a `calc` block.

As the final design comment about Sugawara constructions on a general vector space $V$, it was convenient here to use the literal definition that an $\tilde{\mathbb{F}}$-representation of an $\mathbb{F}$-Lie algebra is a $\mathbb{F}$-Lie algebra homomorphism $\mathfrak{g} \to \mathrm{End}_{\tilde{\mathbb{F}}}(V)$:

```
abbrev LieAlgebra.Representation (F F' : Type*)
    [CommRing F] [CommRing F']
    (g : Type*) [LieRing g] [LieAlgebra F g]
    (V : Type*) [AddCommGroup V] [Module F' V]
    [Module F V] [SMul F F'] [IsScalarTower F F' V]
```

---

[12] The scalar $k \in \mathbb{Z}$ needs to be coerced to an element of the ground field $\mathbb{F}$ by `(k : F)`. The neutral element `1` is that of the ring `V →₁[F] V` of linear endomorphisms of $V$.

[13] The proof states in intermediate steps are less controlled in tactic mode than in `calc` mode, so when something breaks, it is in particular often harder to localize the problem.

```
[SMulCommClass F' F V] :=
  g →₁{F} V →₁[F'] V
```

Note that the first arrow `→₁{F}` in the Lean code excerpt indeed stands for $\mathbb{F}$-Lie algebra homomorphisms, whereas the second one `→₁[F']` is for $\tilde{\mathbb{F}}$-linear maps. Also worth noting is Mathlib's flexible `IsScalarTower` and `SMulCommClass` typeclasses [Wie21]: instead of saying literally that $\tilde{\mathbb{F}}$ extends $\mathbb{F}$, the assumptions here are: $\mathbb{F}$ acts on $\tilde{\mathbb{F}}$, and $V$ is a vector space over both $\mathbb{F}$ and $\tilde{\mathbb{F}}$ in a compatible way.

The convenience of the `g →₁{F} V →₁[F'] V` spelling here is mainly because with the collections $(J_k)_{k \in \mathbb{Z}}$ and $(L_n)_{n \in \mathbb{Z}}$ of operators, we can linearly extend from basis vectors to define Lie algebra homomorphisms $\mathfrak{hei} \to \mathrm{End}_{\mathbb{F}}(V)$ and $\mathfrak{vir} \to \mathrm{End}_{\mathbb{F}}(V)$ once the right commutation relations are available. Generally, we provided the construction of Lie algebra representations from operators with the right commutation relations as `LieAlgebra.representationOfBasis`.

## 12 Verma modules and Fock spaces

Theorem 10.1 is most commonly applied to charged Fock space representations of the Heisenberg algebra, which are a particular case of Verma modules. Verma modules are a fundamental concept in the representation theory of Lie algebras. With a suitable generalization they allow also constructions of fermionic Fock spaces that appear in variants of the Sugawara construction. The present project contains a definition and construction of Verma modules in a generality that allows the fermionic variants, among others. The designs decisions certainly need to be scrutinized before Mathlib pull-requests of this broadly important notion, but the experiment here should be a useful starting point.

The standard mathematical definition is made in the context of a Lie algebra $\mathfrak{g}$ with a **triangular decomposition**

$$\mathfrak{g} = \mathfrak{g}_- \oplus \mathfrak{g}_0 \oplus \mathfrak{g}_+$$

where $\mathfrak{g}_0 \subset \mathfrak{g}$ is an abelian Lie subalgebra, a Cartan subalgebra, and $\mathfrak{g}_\pm \subset \mathfrak{g}$ are Lie subalgebras. Usually one assumes (or proves starting from semisimplicity) that $\mathfrak{g}_\pm$ have root space decompositions to joint eigenspaces of the adjoint actions of elements of $\mathfrak{g}_0$, and that $\mathfrak{g}_\pm$ are nilpotent (in a suitable relaxed sense in for infinite-dimensional Lie algebras).

Let $\mathfrak{g}$ be a Lie algebra with a triangular decomposition $\mathfrak{g} = \mathfrak{g}_- \oplus \mathfrak{g}_0 \oplus \mathfrak{g}_+$, and let $\lambda \in \mathfrak{g}_0^* = \mathrm{Hom}_{\mathbb{F}}(\mathfrak{g}_0, \mathbb{F})$ be a linear functional on the Cartan subalgebra. In this setup one defines weight spaces, highest weight vectors, highest weight modules, and Verma modules as follows. Here we use the language of modules over universal enveloping algebras instead of representations of Lie algebras for notational convenience and because of the chosen generalization below.

**Definition 12.1.** Let $V$ be a left $\mathcal{U}(\mathfrak{g})$-module. A vector $v \in V$ is a **weight vector** of weight $\lambda$ if

$$Hv = \lambda(H)\, v \qquad \text{for all } H \in \mathfrak{g}_0. \tag{20}$$

and a **highest weight vector** of highest weight $\lambda$ if in addition

$$Xv = 0 \qquad \text{for all } X \in \mathfrak{g}_+. \tag{21}$$

The module $V$ is a **highest weight module** with highest weight $\lambda$ if a highest weight vector $v \in V$ of highest weight $\lambda$ exists that generates the entire module, i.e., $V = \mathcal{U}(\mathfrak{g})\, v$.

Verma modules are particular highest weight modules with the following universal property.

**Definition 12.2 (↗).** A highest weight module $M$ over $\mathcal{U}(\mathfrak{g})$ with highest weight $\lambda$ and corresponding highest weight vector $\mathbb{v} \in M$ is said to be a **Verma module** with highest weight $\lambda$ if the following universal property is satisfied: for any $\mathcal{U}(\mathfrak{g})$-module $V$ with a highest weight vector $v$ of highest weight $\lambda$, there exists a unique $\mathcal{U}(\mathfrak{g})$-module homomorphism $f \colon M \to V$ such that $f(\mathbb{v}) = v$.

As usual, the uniqueness up to isomorphism of Verma modules with a given highest weight $\lambda$ follows from this universal property. Existence can be shown in different ways. It would be standard in mathematical literature to "define" (construct) the Verma module as the induced representation of $\mathfrak{g}$ obtained from the one-dimensional representation of the Lie subalgebra $\mathfrak{g}_0 \oplus \mathfrak{g}_+ \subset \mathfrak{g}$ where $\mathfrak{g}_+$ acts as zero and $\mathfrak{g}_0$ acts via the linear functional $\lambda$ (the existence of this one-dimensional representation does require mild assumptions on $\mathfrak{g}_+$). The following alternative construction via the universal enveloping algebra reflects our chosen formalization.

**Lemma 12.3 (↗).** For $\mathfrak{g} = \mathfrak{g}_- \oplus \mathfrak{g}_0 \oplus \mathfrak{g}_+$ and $\lambda \in \mathfrak{g}_0^*$ as above, let $I_\lambda \subset \mathcal{U}(\mathfrak{g})$ be the left ideal generated by the elements of the form $H - \lambda(H)\, 1$ for $H \in \mathfrak{g}_0$ and $X$ for $X \in \mathfrak{g}_+$. View $\mathcal{U}(\mathfrak{g})$ as a left module over itself, so that $I_\lambda$ becomes a submodule. Then the module $M_\lambda = \mathcal{U}(\mathfrak{g})/I_\lambda$ is a Verma module ("the Verma module") with highest weight $\lambda$, with $\mathbb{v}_\lambda := 1 + I_\lambda \in \mathcal{U}(\mathfrak{g})/I_\lambda$ as its highest weight vector.

Two particular cases of Verma modules of Lie algebras are relevant for the main theme of the project.

**Definition 12.4 (↗).** The Virasoro algebra has a triangular decomposition $\mathfrak{vir} = \mathfrak{vir}_- \oplus \mathfrak{vir}_0 \oplus \mathfrak{vir}_+$ where

$$\mathfrak{vir}_\pm = \mathrm{span}_{\mathbb{F}} \left\{ \mathcal{L}_n \mid \pm n > 0 \right\}, \quad \mathfrak{vir}_0 = \mathrm{span}_{\mathbb{F}} \left\{ \mathcal{L}_0, C \right\}.$$

Let $c, h \in \mathbb{F}$. The **Virasoro Verma module** $\mathcal{V}_{c,h}$ of **central charge** $c$ and **conformal weight** $h$ is the Verma module associated with the above triangular decomposition and the highest weight $\lambda \in \mathfrak{vir}_0^*$ defined by $\lambda(C) = c$ and $\lambda(\mathcal{L}_0) = h$. We denote (a choice of) its highest weight vector by $\mathbb{v}_{c,h}$.

**Definition 12.5 (↗).** The Heisenberg algebra has a triangular decomposition $\mathfrak{hei} = \mathfrak{hei}_- \oplus \mathfrak{hei}_0 \oplus \mathfrak{hei}_+$ where

$$\mathfrak{hei}_\pm = \mathrm{span}_{\mathbb{F}} \left\{ \mathcal{J}_k \mid \pm k > 0 \right\}, \quad \mathfrak{hei}_0 = \mathrm{span}_{\mathbb{F}} \left\{ \mathcal{J}_0, \mathcal{K} \right\}.$$

Let $\alpha \in \mathbb{F}$. The **charged Fock space** $\mathcal{F}_\alpha$ of **charge** $\alpha$ is the Verma module for $\mathfrak{hei}$ associated with the above triangular decomposition and the highest weight $\lambda \in \mathfrak{hei}_0^*$ defined by $\lambda(\mathcal{K}) = 1$ and $\lambda(\mathcal{J}_0) = \alpha$. We denote (a choice of) its highest weight vector by $\mathbb{v}_\alpha$ and call this the **vacuum** vector.

From the Poincaré–Birkhoff–Witt (PBW) theorem, it can be seen that a basis of $\mathcal{V}_{c,h}$ is

$$\left(\mathcal{L}_{-n_m} \cdots \mathcal{L}_{-n_1} \mathbb{v}_{c,h}\right)_{m \in \mathbb{N},\ 1 \le n_1 \le \cdots \le n_m}.$$

Similarly, a basis of $\mathcal{F}_\alpha$ is

$$\left(\mathcal{J}_{-k_m} \cdots \mathcal{J}_{-k_1} \mathbb{v}_\alpha\right)_{m \in \mathbb{N},\ 1 \le k_1 \le \cdots \le k_m}.$$

These statements about explicit bases of Verma modules are not yet formalized, since Mathlib does not currently have the PBW theorem. The PBW theorem was perhaps the only major missing desirable prerequisite of the project. A few compromises in the scope were made to avoid relying on it. If a vector $v \in V$ in a representation of an $\mathbb{F}$-algebra $A$ is an eigenvector with eigenvalue $r \in \mathbb{F}$ of a central element $z \in A$, then $z$ acts as scalar multiplication by $r$ on the entire subrepresentation $Av \subset V$ generated by $v$. Central elements of Lie algebras are central in their universal enveloping algebras, so $C$ acts as $c$ id on $\mathcal{V}_{c,h}$, and $\mathcal{K}$ acts as id on $\mathcal{F}_\alpha$.

A representation $V$ of $\mathfrak{hei}$ where $\mathcal{K}$ acts as $\mathrm{id}_V$ satisfies the hypothesis (H), when $\mathsf{J}_k$ is taken to be the action of $\mathcal{J}_k$ on $V$. It is easy to check that the charged Fock space $\mathcal{F}_\alpha$ satisfies also the hypothesis (T), so Theorem 10.1 applies and gives:

**Corollary 12.6** (). *The charged Fock space $\mathcal{F}_\alpha$ carries a representation of the Virasoro algebra $\mathfrak{vir}$ with central charge $c = 1$ so that the action of $\mathcal{L}_n$ is the operator $\mathsf{L}_n$ given by (19).*

An additional calculation, starting from the definitions of $\mathsf{L}_n$ and $\mathbb{v}_\alpha$, reveals that the vacuum vector satisfies $\mathsf{L}_0 \mathbb{v}_\alpha = \frac{1}{2}\alpha^2 \mathbb{v}_\alpha$ and $\mathsf{L}_n \mathbb{v}_\alpha = 0$ for $n > 0$. From the universal property of Verma modules we therefore get:

**Theorem 12.7** (). *There exists a unique homomorphism of $\mathcal{U}(\mathfrak{vir})$-modules $\mathcal{V}_{1, \frac{\alpha^2}{2}} \to \mathcal{F}_\alpha$ such that $\mathbb{v}_{1, \frac{\alpha^2}{2}} \mapsto \mathbb{v}_\alpha$.*

## 13 Formalizing Verma modules

Note that our chosen Definition 12.1 does not include a requirement for a highest weight vector to be nonzero, and correspondingly a highest weight module could be the zero module. Indeed, without some mild extra assumptions (such as nilpotency of $\mathfrak{g}_+$), there might not exist any nonzero highest weight module. Relatedly, in the construction of Lemma 12.3, some further assumptions would be needed to guarantee that $I_\lambda \subset \mathcal{U}(\mathfrak{g})$ is a *proper* ideal.

Let us then make some observations about the definitions, without claiming that they represent good taste in ordinary mathematics. Definitions 12.1 and 12.2 and the construction of Lemma 12.3 admit a generalization to the following setup of very minimal assumptions. In a module $V$ over an associative algebra $A$ with a subset $A' \subset A$ (not necessarily a subalgebra or even a vector subspace) and a (not necessarily linear) function $\eta \colon A' \to \mathbb{F}$, a vector $v \in V$ satisfying $av = \eta(a)v$ for all $a \in A'$ could still function as a "highest weight vector" of "highest weight" $\eta$. This $\eta$ is meant to encode both the joint eigenvalues (20) of the elements of $\mathfrak{g}_0$ and the annihilation property (21) for the elements of $\mathfrak{g}_+$. The chosen setup in our formalization is essentially that, although instead of a subset $A' \subset A$ and a function $\eta \colon A' \to \mathbb{F}$, it is more convenient to use an indexed family $(a_i, \eta_i)_{i \in I}$ of pairs $a_i \in A$ and $\eta_i \in \mathbb{F}$:

```
def vermaIdeal (η : ι → A × F) :
    Submodule A A :=
  Submodule.span A (Set.range
    fun i ↦ (η i).1 - algebraMap F A (η i).2)

def VermaModule (η : ι → A × F) :=
  A ⧸ vermaIdeal η

def VermaModule.hwVec (η : ι → A × F) :
    VermaModule η :=
  Submodule.Quotient.mk 1
```

In formalization, the essential part of the universal property of Verma modules is the construction of a map from a Verma module to another module with a highest weight vector:

```
def VermaModule.universalMap
    (η : ι → A × F) {u : M}
    (u_hwv : ∀ i,
      (η i).1 • u = algebraMap F A (η i).2 • u) :
    VermaModule η →ₗ[A] M := ...

lemma VermaModule.universalMap_hwVec
    (η : ι → A × F) {u : M}
    (u_hwv : ∀ i,
      (η i).1 • u = algebraMap F A (η i).2 • u) :
    universalMap η u_hwv (hwVec η) = u
```

Again we actually separately formulate a *characteristic predicate* `IsVermaModule`, prove that our construction satisfies this predicate and prove uniqueness up to isomorphism of any two modules satisfying this:

```
noncomputable def IsVermaModule.equiv_of_isVermaModule
    (M₁ M₂ : Type u) (hwv₁ : M₁) (hwv₂ : M₂)
    [AddCommGroup M₁] [Module A M₁]
    [AddCommGroup M₂] [Module A M₂]
    (h₁ : IsVermaModule η M₁ hwv₁)
    (h₂ : IsVermaModule η M₂ hwv₂) :
  M₁ ≃ₗ[A] M₂
```

Some justification is needed for our choice above to define Verma modules in generality achieved by dropping all principled guardrails. Two main reasons guided this design. First, even if one is interested only in the usual Lie algebra setup, packaging the triangular decomposition and the linear functional into the definition is heavy. Notably, the component $\mathfrak{g}_-$ is unused in the definitions, so to carry around that data is dubious. Second, giving a definition applicable to associative algebras has real use cases beyond Lie algebras. For example there are representations of Clifford algebras and Lie superalgebras which follow exactly the same pattern. Indeed some fermionic variants of the Sugawara construction are commonly applied to *fermionic Fock spaces*, which will be Verma modules in this generalized sense.

The usual Verma modules of highest weight $\lambda \in \mathfrak{g}_0^*$ over a Lie algebra $\mathfrak{g} = \mathfrak{g}_- \oplus \mathfrak{g}_0 \oplus \mathfrak{g}_+$ are obtained from the above general definition via packaging a triangular decomposition and a highest weight into an indexed collection $(X, \eta_X)_{X \in \mathfrak{g}_0 \oplus \mathfrak{g}_+}$, where $\eta_H = \lambda(H)$ for $H \in \mathfrak{g}_0$ and $\eta_E = 0$ for $E \in \mathfrak{g}_+$. For this case, we introduced the type `TriangularDecomposition g` encoding a direct sum decomposition $\mathfrak{g} = \mathfrak{g}_- \oplus \mathfrak{g}_0 \oplus \mathfrak{g}_+$, and the type `TriangularDecomposition.weight` of linear functionals on its Cartan direct summand $\mathfrak{g}_0$. To convert from a weight to the data $(a_i, \eta_i)_i$ for our generalized notion of Verma modules, the packaging is done by:

Verma modules over Lie algebras then become a special case of the general `VermaModule`. For the further particular case Definitions 12.4 and 12.5 most relevant in the project, we provide abbreviations. For example the triangular decomposition in the former is formalized as `virasoroTri`, and the functional $\lambda \in \mathfrak{vir}_0^*$ is constructed from the parameters $c, h$ in `VirasoroAlgebra.hw`. A Virasoro Verma module `VirasoroVerma c h` then only requires the two parameters `c h : F`, allowing the user to omit specifying the triangular decomposition and constructing the functional.

The formalized version of Theorem 12.7 is a definition of an $\mathcal{U}(\mathfrak{vir})$-module homomorphism $\mathcal{V}_{1, \frac{\alpha^2}{2}} \to \mathcal{F}_\alpha$ together with its property:

Contrary to Section 11, here the most convenient spelling of Lie algebra representations was as modules over the universal enveloping algebra. We provided constructions `UniversalEnvelopingAlgebra.representation` and `Representation.moduleUniversalEnvelopingAlgebra` for converting from the $\mathcal{U}(\mathfrak{g})$-module spelling to the $\mathfrak{g} \to \mathrm{End}_{\tilde{\mathbb{F}}}(V)$ spelling and back. Mathlib in fact has yet another spelling, which we did not use: `LieModule`. Making principled decisions about which spelling to use when and developing a clean API for the occasionally needed translations warrant further consideration.

## 14 Conclusions and outlook

In this project, the following results related to the Virasoro algebra were formalized in Lean: the calculation of the Witt algebra cohomology in degree two (Thm 6.2), the basic bosonic Sugawara construction (Thm 10.1), and the existence of Virasoro module maps from Virasoro Verma modules to charged Fock spaces (Thm 12.7). To state these theorems, we first defined central extensions of Lie algebras, the Virasoro algebra and the Heisenberg algebra, and Verma modules. Much of this should be pull-requested to Mathlib after polishing.

The project has natural continuations to many directions. Defining Lie algebra cohomology in general and connecting the low-degree API could be modeled after Mathlib's existing group cohomology. A natural first result about Lie algebra cohomology could have been Whitehead's lemma on the vanishing of low degree cohomology of finitely semisimple Lie algebras [CE48]. This was omitted, because the project's focus was infinite-dimensional Lie algebras. Highest weight representation theory of Lie algebras is a classical topic in mathematics [Hum72], certainly wanted in Mathlib. It is likely a major undertaking, starting from the missing Poincaré–Birkhoff–Witt theorem. We hope that our design experiment with Verma modules proves helpful there. For Virasoro algebra representation theory, a significant project could aim at a formal proof of the Feigin–Fuks structure theorem for Verma modules over the Virasoro algebra [FF84, VZ90, IK10], which is foundational to numerous topics in conformal field theory.

The most direct planned continuation is variants of the Sugawara construction other than the basic bosonic case. Apart from their significance in their own right, we view them as an interesting experiment of how Lean formalization lends itself to modifying one completed proof to proofs with similar structure but different details. Formalization has been observed to have advantages over informal mathematics for the purposes of such modifications and refactorings, but proofs of nontrivial results still require effort. The Majorana fermion case in the Ramond and Neveu-Schwarz sectors are the next ones planned. The symplectic fermion case [GK99, Kau00] is another relevant one, and the affine Kac–Moody algebra cases [GO85] and coset constructions [GKO86] are interesting for vertex operator algebra theory.

The main lesson from the project was, once again, that Mathlib's current contents, high quality, and unity enable building interesting mathematics. We are deeply grateful to those who build and maintain this remarkable open mathematics library!

## Acknowledgments

This work was supported by the Research Council of Finland project 346309: Finnish Centre of Excellence in Randomness and Structures, "FiRST". I thank Shinji Koshida, Joonas Vättö, and Thomas Wasserman for useful discussions and comments on the manuscript.

# A  Some proofs

We include detailed proofs, together with links to the formalizations, of Theorems 6.2 and 10.1 as Appendices A.1 and A.2, respectively.

## A.1  Proof of Theorem 6.2

Theorem 6.2 is stated in three parts, (a)–(c), and we break down the proof correspondingly.

For part (a), one must verify that (13) defines a 2-cocycle.

**Lemma A.1** (↗). *The Virasoro cocycle is a 2-cocycle,*
$$\omega_{\mathfrak{vir}} \in Z^2(\mathfrak{witt}, \mathbb{F}).$$

*Proof.* By construction, $\omega_{\mathfrak{vir}} \colon \mathfrak{witt} \times \mathfrak{witt} \to \mathbb{F}$ is bilinear. The antisymmetry of (13) on basis elements of the Witt algebra is clear, so $\omega_{\mathfrak{vir}}$ is antisymmetric. It remains to prove the Leibniz rule, i.e., that for $X, Y, X \in \mathfrak{witt}$, we have
$$\omega_{\mathfrak{vir}}(X, [Y, Z]) = \omega_{\mathfrak{vir}}([X, Y], Z) + \omega_{\mathfrak{vir}}(Y, [X, Z]).$$
This formula is trilinear in $X, Y, Z$, so it suffices to verify it for basis vectors $X = \ell_n$, $Y = \ell_m$, $Z = \ell_k$. We calculate
$$\omega_{\mathfrak{vir}}(\ell_n, [\ell_m, \ell_k]) = \omega_{\mathfrak{vir}}(\ell_n, (m-k)\ell_{m+k})$$
$$= (m-k) \frac{n^3 - n}{12} \delta_{n+m+k, 0}.$$
and
$$\omega_{\mathfrak{vir}}([\ell_n, \ell_m], \ell_k) + \omega_{\mathfrak{vir}}(\ell_m, [\ell_n, \ell_k])$$
$$= \omega_{\mathfrak{vir}}((n-m)\ell_{n+m}, \ell_k) + \omega_{\mathfrak{vir}}(\ell_m, (n-k)\ell_{n+k})$$
$$= (n-m) \frac{(n+m)^3 - (n+m)}{12} \delta_{n+m+k, 0}$$
$$+ (n-k) \frac{m^3 - m}{12} \delta_{n+m+k, 0}.$$
Both of the above results are nonzero only if $k = -(n+m)$, in which case $m - k = 2m + n$ and $n - k = 2n + m$, so it suffices to note that
$$(2m + n)(n^3 - n)$$
$$= (n - m)\big((n+m)^3 - (n+m)\big) + (2n + m)(m^3 - m)$$
to verify the cocycle condition (7) for $\omega_{\mathfrak{vir}}$. □

Part (b) is the nontriviality of the cohomology class of the Virasoro cocycle $\omega_{\mathfrak{vir}}$.

**Lemma A.2** (↗). *The cohomology class $[\omega_{\mathfrak{vir}}] \in H^2(\mathfrak{witt}, \mathbb{F})$ of the Virasoro cocycle is nonzero.*

*Proof.* Assume, by way of contradiction, that $\omega_{\mathfrak{vir}} \in B^2(\mathfrak{witt}, \mathbb{F})$, i.e., that $\omega_{\mathfrak{vir}} = \partial \beta$ for some $\beta \in C^1(\mathfrak{witt}, \mathbb{F})$. Then, in particular, for every $n \in \mathbb{Z}$ we would have
$$\omega_{\mathfrak{vir}}(\ell_n, \ell_{-n}) = \beta\big([\ell_n, \ell_{-n}]\big) = 2n \beta(\ell_0).$$
By the definition (13) of $\omega_{\mathfrak{vir}}$, this would imply
$$\frac{n^3 - n}{12} = 2n \beta(\ell_0)$$

for all $n \in \mathbb{Z}$. Considering for example $n = 3$ and $n = 6$, we then get
$$2 = 6 \beta(\ell_0) \qquad \text{and} \qquad \frac{35}{2} = 12 \beta(\ell_0),$$
which obviously yield a contradiction. □

The remaining task of showing that any 2-cocycle is cohomologous to a scalar multiple of the Virasoro cocycle, first make some auxiliary are observations.

**Lemma A.3** (↗). *For any Witt algebra 2-cocycle $\omega \in Z^2(\mathfrak{witt}, \mathbb{F})$ with coefficients in $\mathbb{F}$, and for any $n, m, k \in \mathbb{Z}$, we have*
$$(m - k) \omega(\ell_n, \ell_{m+k}) + (k - n) \omega(\ell_m, \ell_{n+k})$$
$$+ (n - m) \omega(\ell_k, \ell_{n+m}) = 0.$$

*Proof.* Direct from Definition 4.1, (7), and (12). □

**Lemma A.4** (↗). *Let $\omega \in Z^2(\mathfrak{witt}, \mathbb{F})$ be a Witt algebra 2-cocycle such that $\omega(\ell_0, \ell_n) = 0$ for all $n \in \mathbb{Z}$. Then for any $n, m \in \mathbb{Z}$ with $n + m \neq 0$, we have*
$$\omega(\ell_n, \ell_m) = 0.$$

*Proof.* Apply Lemma A.3 with $k = 0$. The last term vanishes, and by skew-symmetry of $\omega$, the first two terms simplify to yield
$$(m + n) \omega(\ell_n, \ell_m) = 0,$$
which, assuming $n + m \neq 0$, yields the asserted equation $\omega(\ell_n, \ell_m) = 0$. □

We are then ready to tackle part (c), too.

**Lemma A.5** (↗). *For any 2-cocycle $\omega \in Z^2(\mathfrak{witt}, \mathbb{F})$, there exists a coboundary $\partial \beta$ with $\beta \in C^1(\mathfrak{witt}, \mathbb{F})$ such that*
$$\omega + \partial \beta = r \cdot \omega_{\mathfrak{vir}}$$
*for some scalar $r \in \mathbb{F}$.*

*Proof.* Let $\omega \in Z^2(\mathfrak{witt}, \mathbb{F})$ be a Witt algebra 2-cocycle. Define a Witt algebra 1-cocycle $\beta \in C^1(\mathfrak{witt}, \mathbb{F})$ by linear extension of
$$\beta(\ell_n) = \begin{cases} -\frac{1}{2} \omega(\ell_1, \ell_{-1}) & \text{if } n = 0 \\ \frac{1}{n} \omega(\ell_0, \ell_n) & \text{if } n \neq 0. \end{cases}$$
For any $n \neq 0$, we calculate
$$(\omega + \partial \beta)(\ell_0, \ell_n) = \omega(\ell_0, \ell_n) + \beta([\ell_0, \ell_n])$$
$$= \omega(\ell_0, \ell_n) - n \beta(\ell_n)$$
$$= \omega(\ell_0, \ell_n) - n \frac{1}{n} \omega(\ell_0, \ell_n) = 0.$$
This property and Lemma A.4 imply that
$$(\omega + \partial \beta)(\ell_0, \ell_n) = 0$$
whenever $n + m \neq 0$.

We will show the asserted equation with
$$r = 2 (\omega + \partial \beta)(\ell_2, \ell_{-2}). \tag{22}$$

By comparison with the Virasoro cocycle $\omega_{\mathfrak{vir}}$ (13), and using skew-symmetry, it remains to show that for any $n \in \mathbb{N}$ we have

$$(\omega + \partial \beta)(\ell_n, \ell_{-n}) = r \frac{n^3 - n}{12}.$$

The case $n = 0$ is a direct consequence of antisymmetry. The case $n = 1$ follows using the definition of $\beta$ and the calculation

$$(\omega + \partial \beta)(\ell_1, \ell_{-1}) = \omega(\ell_1, \ell_{-1}) + \beta([\ell_1, \ell_{-1}])$$
$$= \omega(\ell_1, \ell_{-1}) + 2\beta(\ell_0)$$
$$= \omega(\ell_1, \ell_{-1}) - 2 \frac{1}{2} \omega(\ell_1, \ell_{-1}) = 0.$$

The case $n = 2$ follows directly by the choice (22) of $r$. We prove the equality in the cases $n \geq 3$ by induction on $n$. Assume the equation for smaller values of $n$. Apply Lemma A.3 to $\omega + \partial \beta$ with $m = 1 - n$ and $k = -1$ to get

$$0 = (2 - n)(\omega + \partial \beta)(\ell_n, \ell_{-n}) + (-1 - n)(\omega + \partial \beta)(\ell_{1-n}, \ell_{n-1})$$
$$+ (2n - 1)(\omega + \partial \beta)(\ell_1, \ell_{-1})$$
$$= (2 - n)(\omega + \partial \beta)(\ell_n, \ell_{-n}) - (-1 - n) r \frac{(n-1)^3 - (n-1)}{12}$$
$$= (2 - n)(\omega + \partial \beta)(\ell_n, \ell_{-n}) + \frac{r}{12}(n+1)n(n-1)(n-2).$$

where in the second step we used the induction hypothesis. Since $2 - n \neq 0$, this can be solved for

$$(\omega + \partial \beta)(\ell_n, \ell_{-n}) = -\frac{r}{12} \frac{(n+1)n(n-1)(n-2)}{2-n} = r \frac{n^3 - n}{12},$$

completing the induction step. □

*Proof of Theorem 6.2.* Parts (a), (b), and (c) of the theorem are given by Lemmas A.1, A.5, and A.2, respectively. □

## A.2 Proof of Theorem 10.1

Let $V$ be a vector space $V$ and $(\mathsf{J}_k)_{k \in \mathbb{Z}}$ a collection of operators on $V$. We assume (H) and (T) from Section 10 as needed. Recall also the normal ordered product notation $:\mathsf{J}_k \mathsf{J}_l: \colon V \to V$ from (18), and the definition (19) of the operators $(\mathsf{L}_n)_{n \in \mathbb{Z}}$ in terms of these. Informally speaking, $\mathsf{L}_n$ is $\frac{1}{2} \sum_{k \in \mathbb{Z}} :\mathsf{J}_{n-k} \mathsf{J}_k:$, but this sum contains infinitely many nonzero terms and we do not introduce a topology on $\mathrm{Hom}_{\mathbb{F}}(V, V)$. Instead, acting on any given $v \in V$ we obtain a sum (19) with only finitely many nonzero terms by the assumption (T).

The following alternative form of the normal ordering (18) is a direct consequence of (H), but it makes some rearrangements in calculations nicer:

**Lemma A.6** (↗). *Suppose that $(\mathsf{J}_k)_{k \in \mathbb{Z}}$ satisfy the commutation relations (H). Then for any $k, l \in \mathbb{Z}$ we have*

$$:\mathsf{J}_k \mathsf{J}_l: \; = \; :\mathsf{J}_l \mathsf{J}_k: \; = \; \begin{cases} \mathsf{J}_k \circ \mathsf{J}_l & \text{if } k < 0 \\ \mathsf{J}_l \circ \mathsf{J}_k & \text{if } k \geq 0. \end{cases}$$

Assuming (T), we have the following guarantee that the sum (19) defining $\mathsf{L}_n$ has only finitely many nonzero terms when acting on any given vector is:

**Lemma A.7** (↗). *Suppose that $(\mathsf{J}_k)_{k \in \mathbb{Z}}$ satisfy the local truncation condition (T). Then for any $n \in \mathbb{Z}$ and any $v \in V$, there are only finitely many $k \in \mathbb{Z}$ such that $:\mathsf{J}_{n-k} \mathsf{J}_k: v \neq 0$.*

*Proof.* Fixing $v \in V$, the local truncation condition (T) gives the existence of an $N$ such that $\mathsf{J}_k v = 0$ for $k \geq N$. By the definition (18) of normal ordering, we then have $:\mathsf{J}_k \mathsf{J}_l: v = 0$ when $\max\{k, l\} \geq N$. For a fixed $n \in \mathbb{Z}$, it follows that $:\mathsf{J}_{n-k} \mathsf{J}_k: v = 0$ unless $n - N < k < N$. □

For calculations, we need to observe that the formally infinite sum (19) can be pulled out of commutators.

**Lemma A.8** (↗). *Suppose that $(\mathsf{J}_k)_{k \in \mathbb{Z}}$ satisfy the local truncation condition (T), and let $\mathsf{A} \colon V \to V$ be a linear operator. Then for any $n \in \mathbb{Z}$, the action of the commutator $[\mathsf{L}_n, \mathsf{A}] = \mathsf{L}_n \circ \mathsf{A} - \mathsf{A} \circ \mathsf{L}_n$ on any $v \in V$ is given by the series*

$$[\mathsf{L}_n, \mathsf{A}] v = \frac{1}{2} \sum_{k \in \mathbb{Z}} [:\mathsf{J}_{n-k} \mathsf{J}_k:, \mathsf{A}] v$$

*where only finitely many of the terms are nonzero.*

*Proof.* Write

$$[\mathsf{L}_n, A] v = \mathsf{L}_n A v - A \mathsf{L}_n v$$
$$= \frac{1}{2} \sum_{k \in \mathbb{Z}} :\mathsf{J}_{n-k} \mathsf{J}_k: A v - \frac{1}{2} A \sum_{k \in \mathbb{Z}} :\mathsf{J}_{n-k} \mathsf{J}_k: v.$$

By Lemma A.7, only finitely many of the terms in both sums are nonzero. In particular, in the second term, the linear operator A can be moved inside the sum, to act on the individual terms. The resulting sums with finitely many nonzero terms may be combined, and they become of the asserted form of sum of commutators. □

Our first lemma with actual content is the following: it is one formulation of the physics statement that the current is a primary field of conformal weight 1.

**Lemma A.9** (↗). *Suppose that $(\mathsf{J}_k)_{k \in \mathbb{Z}}$ satisfy the commutation relations (H) and the local truncation condition (T). Then for any $n \in \mathbb{Z}$ and $k \in \mathbb{Z}$, we have*

$$[\mathsf{L}_n, \mathsf{J}_k] \; = \; -k \, \mathsf{J}_{n+k}.$$

*Proof.* Fix $v \in V$. Use Lemma A.8 to write

$$[\mathsf{L}_n, \mathsf{J}_k] v = \frac{1}{2} \sum_{l \in \mathbb{Z}} [:\mathsf{J}_{n-l} \mathsf{J}_l:, \mathsf{J}_k] v.$$

Recall that $:\mathsf{J}_{n-l} \mathsf{J}_l:$ is in any case a product of two factors, only the order depends on the indices $l$ and $n$. We therefore get to use the formula $[AB, C] = A[B, C] + [A, C]B$ for commutators of a product. For example if $n - l \leq l$, we have

$$[:J_{n-l} J_l:, J_k] = [J_{n-l} J_l, J_k]$$
$$= J_{n-l} [J_l, J_k] + [J_{n-l}, J_k] J_l$$
$$= l \, \delta_{l+k,0} J_{n-l} + (n - l) \, \delta_{n-l+k,0} J_l,$$

where we used (H) in the last step. The case $n - l > l$ is similar, and the result is the same. Inserting in the original calculation, we find

$$[L_n, J_k] v = \frac{1}{2} \sum_{l \in \mathbb{Z}} \Big( l \, \delta_{l+k,0} J_{n-l} + (n - l) \, \delta_{n-l+k,0} J_l \Big) v$$
$$= \frac{1}{2} \Big( - k \, J_{n+k} v + (n - (n + k)) \, J_{n+k} v \Big)$$
$$= - k \, J_{n+k} v,$$

completing the proof. □

The next step towards the calculation of $[L_n, L_m]$ is to calculate an individual term $[L_n, :J_k J_{m-k}:]$.

**Lemma A.10** (↗). *Suppose that $(J_k)_{k \in \mathbb{Z}}$ satisfy the commutation relations (H) and the local truncation condition (T). Then for any $n \in \mathbb{Z}$ and $k, m \in \mathbb{Z}$, we have*

$$[L_n, :J_{m-k} J_k:]$$
$$= - k :J_{m-k} J_{n+k}: - (m - k) :J_{n+m-k} J_k:$$
$$+ k (n + k) \, \delta_{n+m,0} \Big( \mathbb{I}_{0 \leq k < -n} - \mathbb{I}_{-n \leq k < 0} \Big) \mathrm{id}_V.$$

*where $\mathbb{I}_{\text{condition}}$ is defined as 1 if the condition is true and 0 otherwise.*

*Proof.* Note again that $:J_{m-k} J_k:$ is a product of two factors, only the order depends on the indices $m$ and $k$. We can therefore use the formula $[A, BC] = B[A, C] + [A, B]C$ for the commutator with a product. Here it is convenient to use Lemma A.6 for the orders of the products.

If $k \geq 0$, we have

$$[L_n, :J_{m-k} J_k:] = [L_n, J_{m-k} J_k]$$
$$= J_{m-k} [L_n, J_k] + [L_n, J_{m-k}] J_k$$
$$= - k \, J_{m-k} J_{n+k} - (m - k) \, J_{n+m-k} J_k.$$

The product in the second term is $J_{n+m-k} J_k = :J_{n+m-k} J_k:$ since $k \geq 0$. The product in the first term differs from a normal ordered product if $n + k < 0$, resulting in a difference which is a commutator,

$$J_{m-k} J_{n+k} - :J_{m-k} J_{n+k}: = \mathbb{I}_{n+k<0} [J_{m-k}, J_{n+k}]$$
$$= \mathbb{I}_{n+k<0} \, \delta_{n+m,0} \, (m - k) \, \mathrm{id}_V$$
$$= - (n + k) \, \delta_{n+m,0} \, \mathbb{I}_{0 \leq k < -n} \, \mathrm{id}_V.$$

Combining the observations, in the case $k \geq 0$ we found

$$[L_n, :J_{m-k} J_k:] = - k :J_{m-k} J_{n+k}: - (m - k) :J_{n+m-k} J_k:$$
$$+ k(n + k) \, \delta_{n+m,0} \, \mathbb{I}_{0 \leq k < -n} \, \mathrm{id}_V.$$

The case $k < 0$ is similar, and the result is

$$[L_n, :J_{m-k} J_k:] = - k :J_{m-k} J_{n+k}: - (m - k) :J_{n+m-k} J_k:$$
$$- k(n + k) \, \delta_{n+m,0} \, \mathbb{I}_{-n \leq k < 0} \, \mathrm{id}_V.$$

In both cases, $k \geq 0$ and $k < 0$, the result matches the assertion. □

Before the final calculation, we record one more auxiliary result.

**Lemma A.11** (↗). *For any $n \in \mathbb{N}$, we have*

$$\sum_{\substack{l \in \mathbb{Z} \\ 0 \leq l < n}} (n - l)l = \frac{n^3 - n}{6}.$$

*Proof.* This is a straightforward induction on $n$. □

The following final calculation completes the proof of Theorem 10.1.

**Lemma A.12** (↗). *Suppose that $(J_k)_{k \in \mathbb{Z}}$ satisfy the commutation relations (H) and the local truncation condition (T). Then for any $n, m \in \mathbb{Z}$, we have*

$$[L_n, L_m] = (n - m) L_{n+m} + \delta_{n+m,0} \frac{n^3 - n}{12} \mathrm{id}_V.$$

*Proof.* Fix $v \in V$. Use Lemmas A.8 and A.10 to write

$$[L_n, L_m] v$$
$$= \frac{1}{2} \sum_{k \in \mathbb{Z}} [L_n, :J_{m-k} J_k:] v$$
$$= \frac{1}{2} \sum_{k \in \mathbb{Z}} \Big( - k :J_{m-k} J_{n+k}: - (m - k) :J_{n+m-k} J_k: \Big) v$$
$$+ \frac{1}{2} \delta_{n+m,0} \sum_{k \in \mathbb{Z}} k (n + k) \Big( \mathbb{I}_{0 \leq k < -n} - \mathbb{I}_{-n \leq k < 0} \Big) v.$$

We claim that the first sum is equal to $(n - m) L_{n+m} v$. Indeed we can split the sum in two, and change variables in the first to $l = n + k$:

$$\sum_{k \in \mathbb{Z}} \Big( - k :J_{m-k} J_{n+k}: v - (m - k) :J_{n+m-k} J_k: v \Big)$$
$$= \sum_{l \in \mathbb{Z}} (n - l) :J_{n+m-l} J_l: v - \sum_{k \in \mathbb{Z}} (m - k) :J_{n+m-k} J_k: v$$
$$= (n - m) \sum_{k \in \mathbb{Z}} :J_{n+k-k} J_k: v \;=\; (n - m) L_{n+m} v.$$

Therefore we have

$$[L_n, L_m] - (n - m) L_{n+m}$$
$$= \frac{1}{2} \delta_{n+m,0} \sum_{k \in \mathbb{Z}} k (n + k) \Big( \mathbb{I}_{0 \leq k < -n} - \mathbb{I}_{-n \leq k < 0} \Big) \mathrm{id}_V,$$

and it suffices to show that

$$\frac{1}{2} \sum_{k \in \mathbb{Z}} k (n + k) \Big( \mathbb{I}_{0 \leq k < -n} - \mathbb{I}_{-n \leq k < 0} \Big) \;=\; \frac{n^3 - n}{12}.$$

Considering cases $n \geq 0$ and $n \leq 0$ separately, this equality reduces to Lemma A.11 in both cases. □